\documentclass[12pt,thmsa]{article}
\usepackage{amsmath}
\usepackage{graphicx}
\usepackage{amsfonts}
\usepackage{amssymb}
\usepackage{mathrsfs}
\newtheorem{theorem}{Theorem}[section]
\newtheorem{lemma}[theorem]{Lemma}
\newtheorem{corollary}[theorem]{Corollary}

\begin{document}

\begin{center}
{\Large \textbf{Dimension reduction in spatial regression with kernel SAVE method}}

\bigskip

M\`etolidji Moquilas Raymond AFFOSSOGBE\textsuperscript{a} , Guy Martial  NKIET\textsuperscript{b}  and Carlos OGOUYANDJOU\textsuperscript{a}

\bigskip

\textsuperscript{a}Institut de Math\'ematiques et de Sciences Physiques, Porto Novo, B\'enin.
\textsuperscript{b}Universit\'{e} des Sciences et Techniques de Masuku,  Franceville, Gabon.

\bigskip

E-mail adresses : metolodji.affossogbe@imsp-uac.org,  guymartial.nkiet@mathsinfo.univ-masuku.com,  ogouyandjou@imsp-uac.org.

\bigskip
\end{center}

\noindent\textbf{Abstract.}We consider   the smoothed version of sliced average variance estimation (SAVE) dimension reduction method for dealing with spatially dependent data that are observations  of a strongly mixing random field. We propose   kernel estimators for   the interest matrix and the effective dimension reduction (EDR) space, and show their  consistency. 

\bigskip

\noindent\textbf{AMS 1991 subject classifications: }62G05, 62G20.

\noindent\textbf{Key words:} SAVE; kernel estimator; spatial data; consistency.

\section{Introduction}
\label{Intro}
\noindent Let us consider the semiparametric regression model introduced by Li \cite{li} and defined as
\begin{equation}\label{model}
Y=g(\beta_{1}^{T}X, \beta_{2}^{T}X,\cdots, \beta_{N}^{T}X,\varepsilon), 
\end{equation}
where $Y$ (resp. $X$) is a random variable with values in $\mathbb{R}$  (resp.  $\mathbb{R}^d$, $d\geq 2$), $N$ is an integer such that $N<d$, the parameters $\beta_1,\beta_2,\cdots,\beta_N$ are $d$-dimensional  linearly independent vectors, $\varepsilon   $ is a random variable that is independent of $X$, and $g$ is an arbitrary unkown function.  The estimation of the space spanned by the  $\beta_k$'s, called  the effective dimension  reduction (EDR) space,  is a crucial issue for achieving reduction dimension. For this problem, Li \cite{li} introduced  the  Sliced Inverse Regression (SIR) method  whereas an alternative  method, called  sliced average variance estimation (SAVE), that  is more comprehensive since it  uses first and second moments was proposed in \cite{Cook}. Smoothed versions of these methods, based on kernel estimators, have been proposed later in \cite{ZhuFang} and \cite{ZhuZhu}. Recently, nonparametric statistical  methods   have evolved  with the existence   of spatially dependent  data. So,   kernel  nonparametric estimation of the   spatial regression function have been studied  (\cite{Li}, \cite{Lu}, \cite{Hallin}, \cite{Carbon1},  \cite{Menezes}, \cite{Attouch}). For dimension reduction in  spatial context,  Loubes and Yao \cite{Loubes} investigated the  kernel SIR   method under strong  mixing conditions. In this note, we study the case of kernel SAVE, which had never been done before. In Section 2, we introduce a kernel   estimate of SAVE based on spatially dependent observations. Then, assumptions and consistency results are given in Section 3. The proofs of theorems are  postponed in Section 4.  
\section{Kernel estimation of SAVE  based on spatial data}\label{sec2}

\noindent  In all of the paper, we assume that $\mathbb{E}\left(\Vert X\Vert^2\right)<+\infty$, where $\Vert\cdot\Vert$ is the usual Euclidean norm of $\mathbb{R}^d$,  and that the covariance matrix $\Sigma $ of $X$ is invertible. Putting  $Z =\Sigma^{-1/2}(X- \mathbb{E}(X))$ and denoting by $Cov(Z | Y)$ the  conditional covariance matrix  of $Z$ conditionally to $Y$,  it is shown in \cite{Cook} that the  EDR  space is fully obtained from  the spectral analysis of the   matrix 
\begin{equation}\label{matrix}
\Gamma:=\mathbb{E}[(I_d-Cov(Z | Y))^ 2]=\textbf{I}_d-2\mathbb{E}(C(Y))+\mathbb{E}\left[ C(Y)^2\right] ,
\end{equation}
where  $\textbf{I}_d$ is the $d\times d$  identity matrix and $C(Y):=Cov(Z|Y)=R(Y)-r(Y)r(Y)^T$, where $R(Y)=\mathbb{E}(ZZ^T\vert Y)$ and $r(Y)=\mathbb{E}(Z|Y)$. From the variance decomposition  theorem, we have that $\textbf{I}_d=\mathbb{E}(C(Y))+\Psi$,   where $\Psi=Cov(\mathbb{E}(r(Y))=\mathbb{E}[r(Y)r(Y)^T]$. Therefore,  $\Gamma=-\textbf{I}_d+2\Psi+\Lambda$, where  $\Lambda=\mathbb{E}[C(Y)^2]$, and the estimation of $\Gamma$ boils down to that of the matrices $\Psi$ and $\Lambda$. From now on, we assume that $Y$ admits a density such that $f(y)>0$    for all $y\in \mathbb{R}$. Let us consider a  stationary random field $\{W_{\textbf{i}},\,\,\textbf{i}\in (\mathbb{N}^*)^L\}$ where  $W_{\textbf{i}}=(Z_{\textbf{i}},Y_{\textbf{i}})$ has the same distribution than  $(Z,Y)$. We suppose that this process is observed on a  region   $\mathcal{I}_{\textbf{n}}=\{{\textbf{i}}=(i_1,i_2,\cdots,i_L)\in \mathbb{Z}^L, 1\leq i_k \leq n_k, k=1,2,\cdots L\}$, where ${\textbf{n}}=(n_1,\cdots,n_L)\in (\mathbb{N}^*)^L$. We put   $\widehat{\textbf{n}}=n_1\times n_2 \cdots\times n_L$ and write $\textbf{n}\rightarrow +\infty$ if $\min\{n_i,i=1,2,\cdots,L\}\rightarrow+\infty$. For defining our estimators, we consider a sequence $(b_{\textbf{n}})$  of strictly positive real  numbers converging to zero as  $\textbf{n}\rightarrow +\infty$, and a kernel function $K$  defined on $\mathbb{R}$. An estimator of $f$ is then given by  $\widehat f_{e_{\textbf{n}}}(y)=\max\{e_{\textbf{n}},\widehat f_{\textbf{n}}(y)\}$,  where $(e_{\textbf{n}})$ is a sequence of strictly positive  real numbers such that $\lim\limits_{{\textbf{n}}\rightarrow +\infty }e_{\textbf{n}}=0$, and $$\widehat f_{\textbf{n}}(y)=\dfrac{1}{\widehat{{\textbf{n}}}b_{{\textbf{n}}}}\sum\limits_{{\textbf{i}} \in\mathcal{I}_{\textbf{n}}}K\left( \dfrac{y-Y_{\textbf{i}}}{b_{\textbf{n}}}\right).$$
Then, we consider
\begin{equation*}
		\widehat m_{\textbf{n}}(y)=\dfrac{1}{\widehat{{\textbf{n}}}b_{{\textbf{n}}}}\sum\limits_{{\textbf{i}} \in \mathcal{I}_{\textbf{n}}}K\left( \dfrac{y-Y_{\textbf{i}}}{b_{\textbf{n}}}\right)\,Z_{\textbf{i}},\,\,
		\widehat M_{\textbf{n}}(y)=\dfrac{1}{\widehat{{\textbf{n}}}b_{{\textbf{n}}}}\sum\limits_{{\textbf{i}} \in \mathcal{I}_{\textbf{n}}}K\left( \dfrac{y-Y_{\textbf{i}}}{b_{\textbf{n}}}\right)\,Z_{\textbf{i}}Z_{\textbf{i}}^T ,\,\,
\end{equation*}
\begin{equation*}
	\widehat r_{\textbf{n}}(y)=\dfrac{\widehat m_{\textbf{n}}(y)}{ \widehat f_{e_{\textbf{n}}}(y)},\,\, \widehat R_{\textbf{n}}(y)=\dfrac{\widehat M_{\textbf{n}}(y)}{\widehat f_{e_{\textbf{n}}}(y)},
	\end{equation*}
and we take as estimator of $\Gamma$ the random matrix
\begin{equation}\label{gamma}
\widehat{\Gamma}_{\textbf{n}}=-\textbf{I}_d+2\widehat{\Psi}_{\textbf{n}}+\widehat{\Lambda}_{\textbf{n}},
\end{equation}
where
\begin{equation*}
		\widehat{\Psi}_{\textbf{n}}=\dfrac{1}{\widehat{{\textbf{n}}}}\sum\limits_{{\textbf{i}}\in \mathcal{I}_{\textbf{n}} }\widehat r_{\textbf{n}}(Y_{\textbf{i}})\widehat r_{\textbf{n}}(Y_{\textbf{i}})^T-\overline{Z}\,\overline{Z}^T,\,\,
		\widehat{\Lambda}_{\textbf{n}}=\dfrac{1}{\widehat{{\textbf{n}}}}\sum\limits_{{\textbf{i}}\in\mathcal{I}_{\textbf{n}} }\widehat{C}_{\textbf{n}}(Y_{\textbf{i}})^2
\end{equation*}
with $\overline{Z}=\dfrac{1}{\widehat{{\textbf{n}}}}\sum\limits_{{\textbf{i}}\in\mathcal{I}_{\textbf{n}} }Z_{\textbf{i}}$ and  $\widehat{C}_{\textbf{n}}(Y_{\textbf{i}})=\widehat R_{\textbf{n}}(Y_{\textbf{i}})-\widehat r_{\textbf{n}}(Y_{\textbf{i}})\widehat r_{\textbf{n}}(Y_{\textbf{i}})^T$. 
\section{Assumptions and asymptotic results}
In order to establish the asymptotic results, the following assumptions will be considered.

\bigskip

\noindent\textbf{Assumption 3.1.} $\Gamma$ is  a positive-definite matrix.
\bigskip

\noindent\textbf{Assumption 3.2.} The kernel  $K$ is a density function with  compact support, is  of order $k$ (where $k\geq 3$) and satisfies  $\int \vert u\vert^k\, K(u)\,du=1$ and $|K(x)-K(y)|\leq C|x-y|$ for some $C>0$.

\bigskip

\noindent\textbf{Assumption 3.3.}  The functions  $f$, $r$ and $R$ belong to  $C^k(\mathbb{R})$ and  $\sup_{y\in\mathbb{R}}|f^{(k)}(y)|$, $\sup_{y\in \mathbb{R}}\Vert   m^{(k)}(y)\Vert   $ and $\sup_{y\in \mathbb{R}}\Vert   M^{(k)}(y)\Vert   $ are bounded, where $m(y)=f(y)\,r(y)$ and   $M(y)=f(y)\,R(y)$.
\bigskip

\noindent\textbf{Assumption 3.4.}  $\sqrt{\widehat{\textbf{n}}}\,\mathbb{E}\left[\Vert   R(Y)\Vert   ^2\textbf{1}_{\{f(Y)\leq e_{\textbf{n}}\}}\right]=o\left(1\right)  $, $\sqrt{\widehat{\textbf{n}}}\,\mathbb{E}\left[\Vert   r(Y)\Vert   ^4\textbf{1}_{\{f(Y)\leq e_{\textbf{n}}\}}\right]=o\left(1\right)$ and 
\noindent $\sqrt{\widehat{\textbf{n}}}\,\mathbb{E}\left[\Vert   R(Y)\Vert   \times \Vert   r(Y)\Vert   ^2\textbf{1}_{\{f(Y)\leq e_{\textbf{n}}\}}\right]=o\left(1\right)$.

\bigskip

\noindent\textbf{Assumption 3.5.}  $\Vert      Z \Vert\leq D$ , where $D$ is a strictly positive constant.

\bigskip

\noindent\textbf{Assumption 3.6.} The process  $\{W_{\textbf{i}},\,\,\textbf{i}\in (\mathbb{Z}^*)^L\}$ is  strongly mixing, i.e. there exists a function  $\chi $ from $\mathbb{R}_{+}$ to itself   satisfying   $\chi(t)\downarrow 0$ as $t\rightarrow +\infty$, such that  for all   subsets $S$ and $S^\prime$ of $(\mathbb{Z}^*)^L$, 
\[
\alpha( \mathcal{B}(S), \mathcal{B}(S^\prime)):=\sup_{A\in  \mathcal{B}(S),\,B\in \mathcal{B}(S^\prime)}|\mathbb{P}(A \cap B)-\mathbb{P}(A)\mathbb{P}(B)|\leq \chi(\delta(S,S^\prime))
\]
where $\mathcal{B}(S)$   (resp. $\mathcal{B}(S^\prime)$) denotes the  Borel  $\sigma$-fields  generated by $\{W_{\textbf{i}},\,\textbf{i}\in S\}$  (resp. $\{W_{\textbf{i}},\,\textbf{i}\in S^\prime\}$) and $\delta(S,S^\prime)$ denotes the  Euclidean distance between $S$ and $S^\prime$. 

\bigskip

\noindent\textbf{Assumption 3.7.}   $ b_{\textbf{n}} \sim \widehat{\textbf{n}}^{-c_1} $ and $e_{\textbf{n}}\sim \widehat{\textbf{n}}^{-c_2}$, where $c_1$ and $c_2$ are real numbers satisfying  $c_1>0$, $0<c_2<\frac{2k-1}{4(2k+1)}$ and  $\dfrac{c_2}{k}+\dfrac{1}{4k}<c_1<\dfrac{1}{2}-2c_2$.

\bigskip
\noindent Putting  $ \phi_{\textbf{n}}=b_{\textbf{n}}^{k}+\dfrac{1}{b_{\textbf{n}}}\sqrt{\dfrac{\log \widehat{\textbf{n}}}{\widehat{\textbf{n}}}}$, we have:

\bigskip
		\begin{theorem}\label{theo1}Under Assumptions 3.2-3.6,  if    $\chi(t)=O(t^{-\theta})$, $t>0$, $\theta>2L$ and $\widehat{\textbf{n}}\,b_{\textbf{n}}^3(\log \widehat{\textbf{n}})^{-1}\rightarrow 0,$ $\widehat{\textbf{n}}\,b_{\textbf{n}}^{\theta_1}(\log \widehat{\textbf{n}})^{-1}\rightarrow +\infty$ with $\theta_1=\dfrac{4L+\theta}{\theta-2L}$, then  we have:
			\begin{align*}
			\widehat{\Gamma}_{\textbf{n}}-\Gamma&=O_p\left( \frac{1}{\sqrt{\widehat{\textbf{n}}}}\right) +O_p\left( \frac{\phi_{\textbf{n}}}{e_{\textbf{n}}}\right) +O_p\left( \frac{\phi_{\textbf{n}}^2}{e_{\textbf{n}}^2}\right) +O_p\left( \frac{\phi_{\textbf{n}}^3}{e_{\textbf{n}}^3}\right) +O_p\left( \frac{\phi_{\textbf{n}}^4}{e_{\textbf{n}}^4}\right) +O_p\left(b_{\textbf{n}}^k+ \frac{\phi_{\textbf{n}}^2}{e_{\textbf{n}}^2}\right).
			\end{align*}
		\end{theorem} 

\bigskip

		\begin{corollary}\label{theo2}
			Under Assumptions 3.2-3.7,   if    $\chi(t)=O(t^{-\theta})$, $t>0$, $\theta>2L$ and $\widehat{\textbf{n}}\,b_{\textbf{n}}^3(\log \widehat{\textbf{n}})^{-1}\rightarrow 0,$ $\widehat{\textbf{n}}\,b_{\textbf{n}}^{\theta_1}(\log \widehat{\textbf{n}})^{-1}\rightarrow +\infty$ with $\theta_1=\dfrac{4L+\theta}{\theta-2L}$, then we have
			$\widehat{\Gamma}_{\textbf{n}}-\Gamma=O_p\left(\widehat{\textbf{n}}^{-1/2}\right)$.
		\end{corollary}

\bigskip
		\noindent For dealing  with the $\widehat{\beta}_j$'s  we  assume that $\tau_1,\tau_2,\cdots,\tau_N$ are    orthonormal eigenvectors of $\Gamma$  associated   with  eigenvalues $\lambda_1,\cdots, \lambda_N$  respectively, such that  $\lambda_1>\lambda_2>\cdots>\lambda_N>0$. Let  $\widehat{\tau}_1,\widehat{\tau}_2,\cdots,\widehat{\tau}_N$ be   orthonormal eigenvectors of $\widehat{\Gamma}_{\textbf{n}}$  associated   with the eigenvalues $\widehat{\lambda}_1,\cdots, \widehat{\lambda}_N$ respectively,  such that  $\widehat{\lambda}_1>\widehat{\lambda}_2>\cdots>\widehat{\lambda}_N>0$. For $j\in\{1,\cdots,N\}$, we have $\beta_j=\Sigma^{-1/2}\tau_j$ and we put  $\widehat{\beta}_j=\widehat{\Sigma}_{\textbf{n}}^{-1/2}\widehat{\tau}_j$, where $\widehat{\Sigma}_{\textbf{n}}=\dfrac{1}{\widehat{\textbf{n}}}\sum\limits_{\textbf{i}\in \mathcal{I}_{\textbf{n}} }\left(X_{\textbf{i}}-\overline{X} \right)\left( X_{\textbf{i}}-\overline{X}\right)^T$    and $\overline{X}=\dfrac{1}{\widehat{\textbf{n}}}\sum\limits_{\textbf{i}\in \mathcal{I}_{\textbf{n}} }X_{\textbf{i}}$. Then, we have:

\bigskip
		\begin{corollary}\label{theo3}
			Under Assumptions 3.1-3.6,   if    $\chi(t)=O(t^{-\theta})$, $t>0$, $\theta>2L$ and $\widehat{\textbf{n}}\,b_{\textbf{n}}^3(\log \widehat{\textbf{n}})^{-1}\rightarrow 0,$ $\widehat{\textbf{n}}\,b_{\textbf{n}}^{\theta_1}(\log \widehat{\textbf{n}})^{-1}\rightarrow +\infty$ with $\theta_1=\dfrac{4L+\theta}{\theta-2L}$, then we have  for any  $j\in\{1,\cdots,N\}$,
			 $\Vert \widehat{\beta}_{j}-\beta_j\Vert=o_p(1)$.
		\end{corollary}
 \section{Proofs}
\subsection{Preliminary results}
	In this section we will give some lemmas necessary to get the proofs of Theorem 3.1. We put
\[
\alpha(t)=\sup_{\textbf{i},\textbf{j}\in \mathbb{R}^L,\,\Vert\textbf{i}-\textbf{j}\Vert=t}\alpha\left(\sigma(W_{\textbf{i}}),\sigma(W_{\textbf{j}})\right).
\]
	\begin{lemma}{\textbf{(See Lemma 6.3 in Loubes and Yao  (2013))}}\\
		Let $(X_{\textbf{n}},\textbf{n}\in \mathbb{N}^N)$ be a centered stationary $\alpha$-mixing process. If there exist $\delta>0$ such that,
		 $\mathbb{E}\left(\Vert X\Vert^{2+\delta}\right)<+\infty$ and $\sum \alpha(\widehat{\textbf{n}})^{\frac{\delta}{2+\delta}}<+\infty$, then $\dfrac{1}{\widehat{\textbf{n}}}\sum\limits_{\textbf{i}\in \mathcal{I}_{\textbf{n}}}X_{\textbf{i}}=O_p(1/\widehat{\textbf{n}})
		$
	\end{lemma}
	\begin{lemma}{(\textbf{Carbon et  al. (2007)})}
		Consider  the sets $S_1,S_2,\cdots,S_p$, each  containing  $m$ sites and such that, for all $i\neq j,$ and for $1\leq i,j\leq p$, $\delta(S_i,S_j)\geq\delta_0$. Let $V_1,V_2,\cdots,V_p$ be a sequence of  real random variables with values in $[a,b]$ and    mesurable with respect to $\mathcal{B}(S_1),\mathcal{B}(S_2),\cdots ,\mathcal{B}(S_p)$ respectively. There exists a sequence of independent random variables $V_{1}^{*},V_{2}^{*},\cdots,V_{p}^{*}$ such that $V_{l}^{*}$ has the same distibution than  $V_l$ and satisfies :
		\begin{align*}
		\sum\limits_{l=1}^{p}\mathbb{E}\left(|V_l-V_{l}^{*}|\right)\leq 2p(b-a)\psi((p-1)m,m)\,\chi(\delta_0).
		\end{align*}
	\end{lemma}

\noindent Note that if the process is stong mixing, then $\psi\equiv 1$ and we have $\sum\limits_{l=1}^{p}\mathbb{E}\left(|V_l-V_{l}^{*}|\right)\leq 2p(b-a)\chi(\delta_0)$.

\bigskip
	\begin{lemma}{\textbf{(See Lemma 6.5 in Loubes and Yao (2013))}}
		Under Assumptions 3.2, 3.3 and 3.6, if    $\alpha(t)=O(t^{-\theta})$, $t>0$, $\theta>2L$ and $\widehat{\textbf{n}}\,b_{\textbf{n}}^3(\log \widehat{\textbf{n}})^{-1}\rightarrow 0,$ $\widehat{\textbf{n}}\,b_{\textbf{n}}^{\theta_1}(\log \widehat{\textbf{n}})^{-1}\rightarrow +\infty$ with $\theta_1=\dfrac{4L+\theta}{\theta-2L}$, then 
		\begin{align*}
		\sup_{y\in\mathbb{R} }|f_{\textbf{n}}(y)-f(y)|&=O_p\left( b_{\textbf{n}}^{k}+\sqrt{\dfrac{\log \widehat{{\textbf{n}}}}{\widehat{{\textbf{n}}}b_{\textbf{n}}}}\right)=O_p\left(b_{\textbf{n}}^{k}+\dfrac{1}{b_{\textbf{n}}}\sqrt{\dfrac{\log \widehat{\textbf{n}}}{\widehat{\textbf{n}}}} \right),\\
		\sup_{y\in\mathbb{R} }\Vert\widehat{m}_{\textbf{n}}(y)-m(y)\Vert&=O_p\left( b_{\textbf{n}}^{k}+\sqrt{\dfrac{\log \widehat{{\textbf{n}}}}{\widehat{{\textbf{n}}}b_{\textbf{n}}}}\right)=O_p\left(b_{\textbf{n}}^{k}+\dfrac{1}{b_{\textbf{n}}}\sqrt{\dfrac{\log \widehat{\textbf{n}}}{\widehat{\textbf{n}}}} \right).
		\end{align*}
	\end{lemma}
	
	\begin{lemma}
		Under Assumptions 3.2, 3.3 and 3.6, if    $\alpha(t)=O(t^{-\theta})$, $t>0$, $\theta>2L$ and $\widehat{\textbf{n}}\,b_{\textbf{n}}^3(\log \widehat{\textbf{n}})^{-1}\rightarrow 0,$ $\widehat{\textbf{n}}\,b_{\textbf{n}}^{\theta_1}(\log \widehat{\textbf{n}})^{-1}\rightarrow +\infty$ with $\theta_1=\dfrac{4L+\theta}{\theta-2L}$, then 
		\begin{align*}
		\sup_{y\in\mathbb{R} }\Vert\widehat{M}_{\textbf{n}}(y)-M(y)\Vert&=O_p\left(b_{\textbf{n}}^{k}+\dfrac{1}{b_{\textbf{n}}}\sqrt{\dfrac{\log \widehat{\textbf{n}}}{\widehat{\textbf{n}}}} \right)
		\end{align*}
	\end{lemma}
	\noindent\textit{Proof.} Clearly,
	
		\begin{align*}
		\sup_{y\in \mathbb{R}}\Vert\widehat{M}_{\textbf{n}}(y)-M(y)\Vert&\leq \sup_{y\in \mathbb{R}}|\widehat{M}_{\textbf{n}}(y)-\mathbb{E}[\widehat{M}_{\textbf{n}}(y)]\Vert+\sup_{y\in \mathbb{R}}\Vert\mathbb{E}[\widehat{M}_{\textbf{n}}(y)]-M(y)\Vert. 
		\end{align*}
		$\bullet$  Let us first study  $\sup_{y\in \mathbb{R}}\Vert\mathbb{E}[\widehat{M}_{\textbf{n}}(y)]-M(y)\Vert$.
		\begin{eqnarray*}
		\mathbb{E}[\widehat{M}_{\textbf{n}}(y)]-M(y)&=&\dfrac{1}{\widehat{\textbf{n}}b_{\textbf{n}}}\sum\limits_{{\textbf{i}} \in \mathcal{I}_{\textbf{n}}}\mathbb{E}\left[ Z_{\textbf{i}}Z_{\textbf{i}}^TK\left( \dfrac{y-Y_{\textbf{i}}}{b_{\textbf{n}}}\right) \right] -M(y)\\ 
&=&\dfrac{1}{\widehat{\textbf{n}}b_{\textbf{n}}}\sum\limits_{{\textbf{i}} \in\mathcal{I}_{\textbf{n}}}\mathbb{E}\left[ \mathbb{E}[Z_{\textbf{i}}Z_{\textbf{i}}^T|Y_{\textbf{i}}]K\left( \dfrac{y-Y_{\textbf{i}}}{b_{\textbf{n}}}\right) \right]-M(y)\\
& =& \dfrac{1}{\widehat{\textbf{n}}b_{\textbf{n}}}\sum\limits_{{\textbf{i}} \in\mathcal{I}_{\textbf{n}}}\mathbb{E}\left[ R(Y_{\textbf{i}})K\left( \dfrac{y-Y_{\textbf{i}}}{b_{\textbf{n}}}\right) \right] -M(y)\\
		&=&\dfrac{1}{\widehat{\textbf{n}}b_{\textbf{n}}}\sum\limits_{{\textbf{i}} \in \mathcal{I}_{\textbf{n}}}\int R(u)f(u)K\left( \dfrac{y-u}{b_{\textbf{n}}}\right)du-M(y)\\
		&=&\dfrac{1}{b_{\textbf{n}}}\int R(u)f(u)K\left( \dfrac{y-u}{b_{\textbf{n}}}\right)du-M(y)\\
		&=&\dfrac{1}{b_{\textbf{n}}}\int M(u)K\left( \dfrac{y-u}{b_{\textbf{n}}}\right)du-M(y)\\
		&=&\int K\left(v\right)\left[ M(y-vb_{\textbf{n}})-M(y)\right] dv\\
		&=&\int K\left(v\right)\left[\sum\limits_{i=1}^{k-1}\dfrac{(-vb_{\textbf{n}})^i}{i!}M^{(i)}(y)+\dfrac{(-vb_{\textbf{n}})^k}{k!}M^{(k)}(y-\eta v b_{\textbf{n}})\right] dv,
\end{eqnarray*}
where  $0<\eta<1$. Hence
\begin{eqnarray*}
		\Vert\mathbb{E}[\widehat{M}_{\textbf{n}}(y)]-M(y)\Vert&\leq &\int \left| K\left(v\right)\dfrac{(-vb_{\textbf{n}})^k}{k!}\right| \Vert M^{(k)}(y-\eta v b_{\textbf{n}})\Vert  dv\\
&\leq&
		  \dfrac{b_{\textbf{n}}^k}{k!}\sup_{y\in \mathbb{R}}\Vert M^{(k)}(y)\Vert\int |v|^k K\left(v\right)dv 
		= C_1b_{\textbf{n}}^{k}
\end{eqnarray*}
		and, therefore, $\sup_{y\in \mathbb{R}}\Vert\mathbb{E}[\widehat{M}_{\textbf{n}}(y)]-M(y)\Vert=O_p(b_{\textbf{n}}^{k})$.

\bigskip

\noindent $\bullet$  Secondly,  let us study  $\sup_{y\in \mathbb{R}}\Vert\widehat{M}_{\textbf{n}}(y)-\mathbb{E}[\widehat{M}_{\textbf{n}}(y)]\Vert$.

\bigskip

	\noindent Consider a real $\varepsilon>0$ and  a sequence $(a_{\textbf{n}})$  of non-negative real numbers converging to $+\infty$. We have:
	\begin{align*}
	\mathbb{P}\left( \sup_{y\in \mathbb{R}}\Vert	 \widehat{M}_{\textbf{n}}(y)-\mathbb{E}\left[\widehat{M}_{\textbf{n}}(y) \right]\Vert>\varepsilon\right)
	&=\mathbb{P}\left( \sup_{y\in \mathbb{R}}\Vert	\widehat{M}_{\textbf{n}}(y)-\mathbb{E}\left[ \widehat{M}_{\textbf{n}}(y) \right]\Vert>\varepsilon;\Vert ZZ^T\Vert\leq a_{\textbf{n}}\right)\\
	& +\mathbb{P}\left(\sup_{y\in \mathbb{R}}\Vert	 \widehat{M}_{\textbf{n}}(y)-\mathbb{E}\left[ \widehat{M}_{\textbf{n}}(y) \right]\Vert>\varepsilon;\Vert ZZ^T\Vert> a_{\textbf{n}}\right)\\
	&\leq \mathbb{P}\left( \sup_{y\in \mathbb{R}}\Vert	\widehat{M}_{\textbf{n}}(y)-\mathbb{E}\left[ \widehat{M}_{\textbf{n}}(y) \right]\Vert>\varepsilon;\Vert ZZ^T\Vert\leq a_{\textbf{n}}\right)\\
	&+ \mathbb{P}\left(\Vert ZZ^T\Vert> a_{\textbf{n}}\right).
	\end{align*}
	Under Assumtion 3.2, $K$ is bounded by a constant $C_2>0$. Then we have for any $ y\in \mathbb{R}$
	\begin{align*}
	\Vert	\widehat{M}_{\textbf{n}}(y)-\mathbb{E}\left[\widehat{M}_{\textbf{n}}(y) \right]\Vert&=\bigg\Vert\dfrac{1}{\widehat{\textbf{n}}b_{\textbf{n}}}\sum\limits_{\textbf{i}\in \mathcal{I}_{\textbf{n}}}\left[ Z_{\textbf{i}}Z_{\textbf{i}}^TK\left(\dfrac{Y_i-y}{h} \right) -\mathbb{E}[ Z_{\textbf{i}}Z_{\textbf{i}}^TK\left(\dfrac{Y_i-y}{h}\right)] \right] \bigg\Vert\\
	&\leq \dfrac{C_2}{\widehat{\textbf{n}}b_{\textbf{n}}}\sum\limits_{\textbf{i}\in \mathcal{I}_{\textbf{n}}}\{\Vert  Z_{\textbf{i}}Z_{\textbf{i}}^T\Vert +\mathbb{E}[ \Vert  Z_{\textbf{i}}Z_{\textbf{i}}^T \Vert \}.
\end{align*}
Thus
\[
	\sup_{y\in\mathbb{R}}\Vert	\widehat{M}_{\textbf{n}}(y)-\mathbb{E}\left[\widehat{M}_{\textbf{n}}  (y) \right]\Vert\leq \dfrac{C_2}{\widehat{\textbf{n}}b_{\textbf{n}}}\sum\limits_{\textbf{i}\in\mathcal{I}_{\textbf{n}}}\{\Vert  Z_{\textbf{i}}Z_{\textbf{i}}^T\Vert +\mathbb{E}[ \Vert  Z_{\textbf{i}}Z_{\textbf{i}}^T \Vert \}
\]	
	and
	\begin{align*}
	& \mathbb{P}\left( \sup_{y\in \mathbb{R}}\Vert	\widehat{M}_{\textbf{n}}(y)-\mathbb{E}\left[\widehat{M}_{\textbf{n}}(y) \right]\Vert>\varepsilon;\Vert ZZ^T\Vert\leq a_{\textbf{n}}\right)\\
	&\leq \mathbb{P}\left(  \dfrac{C_2}{\widehat{\textbf{n}}b_{\textbf{n}}}\sum\limits_{\textbf{i}\in \mathcal{I}_{\textbf{n}}}\{\Vert  Z_{\textbf{i}}Z_{\textbf{i}}^T\Vert +\mathbb{E}[ \Vert  Z_{\textbf{i}}Z_{\textbf{i}}^T \Vert \}>\varepsilon;\Vert ZZ^T\Vert\leq a_{\textbf{n}}\right)\\
	&\leq \mathbb{P}\left(  \dfrac{C_2}{\widehat{\textbf{n}}b_{\textbf{n}}}\sum\limits_{\textbf{i}\in\mathcal{I}_{\textbf{n}}}\{\Vert  Z_{\textbf{i}}Z_{\textbf{i}}^T\Vert +\mathbb{E}[ \Vert  Z_{\textbf{i}}Z_{\textbf{i}}^T \Vert \}\textbf{1}_{\{\Vert Z_{\textbf{i}}Z_{\textbf{i}}^T\Vert\leq a_{\textbf{n}}\}}>\varepsilon\right)\\
	&=\mathbb{P}(S_{\textbf{n}}> \varepsilon),
\end{align*}
where 
\[
S_{\textbf{n}}= \dfrac{C_2}{\widehat{\textbf{n}}b_{\textbf{n}}}\sum\limits_{\textbf{i}\in\mathcal{I}_{\textbf{n}}}\{\Vert  Z_{\textbf{i}}Z_{\textbf{i}}^T\Vert +\mathbb{E}[ \Vert  Z_{\textbf{i}}Z_{\textbf{i}}^T \Vert] \}\textbf{1}_{\{\Vert Z_{\textbf{i}}Z_{\textbf{i}}^T\Vert\leq a_{\textbf{n}}\}}=\sum\limits_{\textbf{i}\in\mathcal{I}_{\textbf{n}}}\Theta_{\textbf{i,n}}.
	\]
		As in  Dabo-Niang  et  al.  (2014) we suppose that  $n_\ell=2pq_\ell$ for $1 \leq \ell \leq L$ and we will use the spatial block decomposition of  Tran (1990). By this method, the random variables $\Theta_{\textbf{i,n}}$, $\textbf{i}\in\mathcal{I}_{\textbf{n}}$ can be grouped in $2^Lq_1\cdots q_L$ cubic blocks of side $p$, coming from a partition of $\mathcal{I}_\textbf{n}$. For   $\textbf{m}=(m_1,\cdots,m_L)\in [\![0;q_1-1]\!]\times[\![0;q_2-1]\!] \times \cdots \times [\![0;q_L-1]\!]$ and for $l\in [\![1;2^L]\!]$, we put
		\[
		R_{l,\textbf{m}}=A_{1,l,\textbf{m}}\times A_{2,l,\textbf{m}}\times \cdots \times A_{L,l,\textbf{m}},
\]
where 
\[
A_{k,l,\textbf{m}}\in \{[\![2m_kp+1;(2m_k+1)p]\!],[\![(2m_k+1)p+1;2(m_k+1)p]\!]\}
\]
 with $j\in  [\![1;L]\!]$. Since
\begin{eqnarray*}
[\![1,n_j]\!]&=&[\![1,p]\!]\cup[\![p+1,2p]\!]\cup[\![2p+1,3p]\!]\cup \cdots \cup [\![p(2q_j-1)+1,2pq_j]\!]\\
		&=& \bigcup_{m_j=0}^{q_j-1}\left( [\![ 2m_jp+1,(2m_j+1)p]\!]\cup[\![(2m_j+1)p+1;2(m_j+1)p]\!]\right) ,
\end{eqnarray*}
it follows
		\begin{align*}
		\mathcal{I}_{\textbf{n}}&=\prod_{j=1}^{L}[\![1,n_j]\!]=\bigcup_{l=1}^{2^L}\,\,\,\bigcup_{\textbf{m}\in \prod_{j=1}^{L}[\![0;q_j-1]\!] }R_{l,\textbf{m}}
		\end{align*}
For $\textbf{m}\in \prod_{j=1}^{L}[\![0;q_j-1]\!]$, let us consider
		\begin{align*}
		U(1,\textbf{n},\textbf{m})&=\sum\limits_{\substack{i_j=2m_jp+1,\\j=1,\cdots,L}}^{(2m_j+1)p}\Theta_{\textbf{i,n}}\\
		U(2,\textbf{n},\textbf{m})&=\sum\limits_{\substack{i_j=2m_jp+1,\\j=1,\cdots,L-1}}^{(2m_j+1)p}\sum\limits_{i_L=(2m_L+1)p+1}^{2(m_L+1)p}\Theta_{\textbf{i,n}}\\
		U(3,\textbf{n},\textbf{m})&=\sum\limits_{\substack{i_j=2m_jp+1,\\j=1,\cdots,L-2}}^{(2m_j+1)p}\sum\limits_{i_{L-1}=(2m_{L-1}+1)p+1}^{2(m_{L-1}+1)p}\sum\limits_{i_L=2m_Lp+1}^{(2m_L+1)p}\Theta_{\textbf{i,n}}\\
		U(4,\textbf{n},\textbf{m})&=\sum\limits_{\substack{i_j=2m_jp+1,\\j=1,\cdots,L-2}}^{(2m_j+1)p}\sum\limits_{i_{L-1}=(2m_{L-1}+1)p+1}^{2(m_{L-1}+1)p}\sum\limits_{i_{L}=(2m_{L}+1)p+1}^{2(m_L+1)p}\Theta_{\textbf{i,n}}
\end{align*}
and so on, until
\[
		U(2^{L-1},\textbf{n},\textbf{m})=\sum\limits_{\substack{i_j=(2m_j+1)p+1,\\j=1,\cdots,L-1}}^{2(m_j+1)p}\sum\limits_{i_L=2m_Lp+1}^{(2m_L+1)p}\Theta_{\textbf{i,n}}\,\,\textrm{ and }\,\,
		U(2^L,\textbf{n},\textbf{m})=\sum\limits_{\substack{i_j=(2m_j+1)p+1,\\j=1,\cdots,L}}^{2(m_j+1)p}\Theta_{\textbf{i,n}}.
		\]
		For $1\leq q \leq 2^L$, we consider 
		$T(q,\textbf{n})=\sum\limits_{\substack{m_j=0,\\j=1,\cdots,L}}^{q_j-1}U(q,\textbf{n},\textbf{m})$, and we have $S_{\textbf{n}}=\sum\limits_{q=1}^{2^L}T(q,\textbf{n})$. Then,
\begin{align*}
		\mathbb{P}\left( S_{\textbf{n}}> \varepsilon \right) &\leq  \mathbb{P}\left(\sum\limits_{q=1}^{2^L}|T(q,\textbf{n})|\geq \varepsilon \right) 
		\leq \mathbb{P}\left((|T(1,\textbf{n})|> \varepsilon/2^L) \cup\cdots\cup(|T(2^L,\textbf{n})|> \varepsilon/2^L)\right) \\
		&\leq 2^L \mathbb{P}\left(|T(1,\textbf{n})|> \varepsilon/2^L \right)
		\end{align*}
			as  the $T(q,\textbf{n})$'s, $1\leq q \leq 2^L$,  have the same distribution. Let denote  $\widetilde{q}=q_1\times \cdots \times q_L$ and $V_1,V_2,\cdots,V_{\widetilde{q}}$, the $\widetilde{q}$ terms  $U(1,\textbf{n},\textbf{m}),\textbf{m}\in [\![1,q_1]\!]\times \cdots \times [\![1,q_L]\!]$ of the sum $T(1,\textbf{n})$. Then $T(1,\textbf{n})=\sum\limits_{l=1}^{\widetilde{q}}V_l$. Let remark that each $\textbf{m}\in [\![1,q_1]\!]\times \cdots \times[\![1,q_L]\!]$, $U(1,\textbf{n},\textbf{m})$ is mesurable with respect to the sigma algebra spanned by the $\Theta_{\textbf{i,n}}$ where $2m_jp+1\leq i_j \leq (2m_j+1)p, j=1,2,\cdots,L$. The sets of those sites are separated  by a distance at least equal to $p$. Indeed, Let $\textbf{m}\in \prod_{j=1}^{L}[\![0;q_j-1]\!]$ and $\textbf{m}^{'}\in \prod_{j=1}^{L}[\![0;q_j-1]\!]$ such that $\textbf{m}^{'} \neq \textbf{m}$, then there exist $j\in [\![1;L]\!]$ such that $m_j=m^{'}_j+u,\,u\in \mathbb{N}^{*}$. Denoting by   $E_{\textbf{m}}$ and $E_{\textbf{m}^{'}}$ the set of sites associated to $U(1,\textbf{n},\textbf{m})$ and $U(1,\textbf{n},\textbf{m}^{'})$ respectively,  we have 
			\begin{align*}
			E_{\textbf{m}}&=\{\textbf{i}=(i_1,i_2,\cdots,i_L)  /i_j\in [\![2m_jp+1;(2m_j+1)]\!],\,j\in [\![1;L]\!]\}\\
			E_{\textbf{m}^{'}}&=\{\textbf{i}^{'}=(i^{'}_1,i^{'}_2,\cdots,i^{'}_L) /i_j\in [\![2m^{'}_jp+1;(2m^{'}_j+1)]\!],j\in [\![1;L]\!]\}
                               \end{align*}
and, therefore,
\begin{eqnarray*}
			\delta(	E_{\textbf{m}},	E_{\textbf{m}^{'}})&=&\min\{\Vert\textbf{i}-\textbf{i}^{'}\Vert,\textbf{i}\in E_{\textbf{m}} , \textbf{i}^{'}\in  E_{\textbf{m}^{'}}\}
\geq \vert i_j-i_{j}^{'}|\\
&=&|2m_jp+l-2m_{j}^{'}p+l^{'}\vert
			=\vert 2up+l-l^{'}\vert
\end{eqnarray*}
 where $(l,l^{'})\in [\![1;p-1]\!]^2$. Since $l-l^{'}+ pu\geq 0$  and  $up\geq 0$, it follows: $\delta(	E_{\textbf{m}},	E_{\textbf{m}^{'}})\geq pu\geq p$.
In addition for all $l\in[\![1,\widetilde{q}]\!]$,
			\begin{align*}
			V_l&\leq \sum\limits_{\substack{i_j=2m_jp+1,\\j=1,\cdots,L}}^{(2m_j+1)p}\Theta_{\textbf{i,n}}\leq \dfrac{C_2}{\widehat{\textbf{n}}b_{\textbf{n}}}\sum\limits_{\substack{i_j=2m_jp+1,\\j=1,\cdots,L}}^{(2m_j+1)p}\{\Vert  Z_{\textbf{i}}Z_{\textbf{i}}^T\Vert +\mathbb{E}[ \Vert  Z_{\textbf{i}}Z_{\textbf{i}}^T \Vert] \}\textbf{1}_{\{\Vert Z_{\textbf{i}}Z_{\textbf{i}}^T\Vert\leq a_{\textbf{n}}\}}\leq \dfrac{2C_2p^La_{\textbf{n}}}{\widehat{\textbf{n}}b_{\textbf{n}}}\leq \dfrac{C_2a_{\textbf{n}}}{b_{\textbf{n}}}.
			\end{align*}
			Then, from Lemma 2,  there exist i.i.d random variables  
			$V_{1}^{*},\cdots,V_{\widetilde{q}}^{*}$ such that for any $l\in[\![1,\widetilde{q}]\!]$, $V_{l}^{*}$ have the  same distribution than  $V_{l}$,  and
			\begin{align*}
			\sum\limits_{l=1}^{\widetilde{q}}\mathbb{E}(|V_l-V_{l}^{*}|)&\leq 2\widetilde{q}\dfrac{p^La_{\textbf{n}}}{\widehat{\textbf{n}}b_{\textbf{n}}}\chi(p)
			\leq \dfrac{2\widetilde{q}p^{L-\theta}a_{\textbf{n}}}{\hat{\textbf{n}}b_{\textbf{n}}}.
			\end{align*}
On the other hand, 
			\begin{align*}
			\mathbb{P}(S_{\textbf{n}}> \varepsilon)&\leq 2^LP\left( \sum\limits_{l=1}^{\widetilde{q}}V_l>\dfrac{\varepsilon}{2^L}\right) 
			\leq 2^L  \mathbb{P}\left( \sum\limits_{l=1}^{\widetilde{q}}\vert 
			V_l-V_{l}^{*}+V_{l}^{*}\vert> \dfrac{\varepsilon}{2^L}\right) \\
			&\leq  2^L  \mathbb{P}\left( \sum\limits_{l=1}^{\widetilde{q}}\vert V_l-V_{l}^{*}\vert>\dfrac{\varepsilon}{2^{L+1}}\right)+ 2^L  \mathbb{P}\left( \sum\limits_{l=1}^{\widetilde{q}}\vert V_{l}^{*}\vert> \dfrac{\varepsilon}{2^{L+1}}\right).
			\end{align*}
			Then, using the Markov inequality, we obtain
			\begin{align*}
			\mathbb{P}\left( \sum\limits_{l=1}^{\widetilde{q}}\vert V_l-V_{l}^{*}\vert> \dfrac{\varepsilon}{2^{L+1}}\right)&\leq\sum\limits_{1}^{\widetilde{q}}\dfrac{2^{L+1}}{\varepsilon}\mathbb{E}(\vert V_l-V_{l}^{*}\vert)
			\leq \dfrac{2^{L+2}\widetilde{q}p^{L-\theta}a_{\textbf{n}}}{\varepsilon \widehat{\textbf{n}}b_{\textbf{n}}}\leq \dfrac{C_3p^{-\theta}a_{\textbf{n}}}{\varepsilon b_{\textbf{n}}},
			\end{align*}
where $C_3>0$.  Taking $\varepsilon=\varepsilon_{\textbf{n}}=\dfrac{1}{b_{\textbf{n}}}\sqrt{\dfrac{\log \widehat{\textbf{n}}}{\widehat{\textbf{n}}}}$,  $p=p_{\textbf{n}}=\left[ \left( \dfrac{\log \widehat{\textbf{n}}}{\widehat{\textbf{n}}b_{\textbf{n}}}\right)^{-1/2L}\right] \sim \left( \dfrac{\log \widehat{\textbf{n}}}{\widehat{\textbf{n}}b_{\textbf{n}}}\right)^{-1/2L}$ and $a_{\textbf{n}}=(\log \widehat{\textbf{n}})^{1/4}$ leads to  $
			a_{\textbf{n}}b_{\textbf{n}}^3\varepsilon=\left( (\log \widehat{\textbf{n}})^{3/2}b_{\textbf{n}}^4\widehat{\textbf{n}}^{-1}\right) ^{1/2}$. Since, from  the hypotheses of the lemma, we have
			$(\log \widehat{\textbf{n}})^{-1}b_{\textbf{n}}^3\widehat{\textbf{n}}\rightarrow 0 $ as $\textbf{n}\rightarrow +\infty$,   there exist a constant $C_4>0$, such that $b_{\textbf{n}}^3\leq C_4\dfrac{\log \widehat{\textbf{n}} }{\widehat{\textbf{n}}}$ for $\textbf{n}$  large  enough. Therefore, we have
			$a_{\textbf{n}}b_{\textbf{n}}^3\varepsilon\leq\sqrt{C_4}b_{\textbf{n}}^{1/2}\left( \dfrac{\log \widehat{\textbf{n}} }{\widehat{\textbf{n}}^{4/5}}\right) ^{5/4}$ from what we deduce that $a_{\textbf{n}}b_{\textbf{n}}^3\varepsilon\rightarrow 0$  as $\textbf{n}\rightarrow +\infty$. Then, exists a constant $C_5>0$ such that $a_{\textbf{n}}\leq C_5(b_{\textbf{n}}^3\varepsilon)^{-1}$ for  $\textbf{n}$  large enough. Consequently, we obtain the inequality
			\begin{align*}
			\dfrac{C_3p^{-\theta}a_{\textbf{n}}}{\varepsilon b_{\textbf{n}}}
			\leq C_3C_5p^{-\theta}\varepsilon^{-2}b_{\textbf{n}}^{-4}\leq C_6\left( \dfrac{\log \widehat{\textbf{n}} }{\widehat{\textbf{n}}b_{\textbf{n}}}\right) ^{\theta/2L}\left( \dfrac{1}{b_{\textbf{n}}}\sqrt{\dfrac{\log \widehat{\textbf{n}} }{\widehat{\textbf{n}}}}\right)^{-2}b_{\textbf{n}}^{-4}\leq  C_6\left( \widehat{\textbf{n}}(\log \widehat{\textbf{n}})^{-1} b_{\textbf{n}}^{\theta_1}\right) ^{\frac{2L-\theta}{2L}}.
			\end{align*}
			In other side, using  Bernstein inequality for i.i.d. bounded random variable, we obtain
			\begin{align*}
			\mathbb{P}\left( \sum\limits_{l=1}^{\widetilde{q}} V_{l}^{*}> \dfrac{\varepsilon}{2^{L+1}}\right)&\leq \exp\left(-\dfrac{\widehat{\textbf{n}}\varepsilon^2b_{\textbf{n}}^2}{4C_{2}^2a_{\textbf{n}}^2} \right)  \leq \exp(-C_7(\log \widehat{\textbf{n}})^{1/2})
			\end{align*}
and we deduce that  $\mathbb{P}\left( \sum\limits_{l=1}^{\widetilde{q}} V_{l}^{*}> \dfrac{\varepsilon}{2^{L+1}}\right)\rightarrow 0$ as $\textbf{n}\rightarrow +\infty$. 
			In addition,  from Markov inequality, we have
			\begin{align*}
			\mathbb{P}(\Vert ZZ^T \Vert>a_{\textbf{n}})&\leq \dfrac{\mathbb{E}[\Vert ZZ^T\Vert]}{a_{\textbf{n}}}= \dfrac{\mathbb{E}[\Vert ZZ^T\Vert]}{(\log(\widehat{\textbf{n}}))^{1/4}}
\end{align*}
and we conclude that 
			$\mathbb{P}(\Vert ZZ^T \Vert>a_{\textbf{n}})\rightarrow 0$ as $\textbf{n}\rightarrow +\infty$. From all what precede, we can conlude that
			$\sup_{y\in \mathbb{R}}\Vert\mathbb{E}[\widehat{M}_{\textbf{n}}(y)]-\widehat{M}_{\textbf{n}}(y)\Vert=O_p\left(  \dfrac{1}{b_{\textbf{n}}}\sqrt{\dfrac{\log \widehat{\textbf{n}}}{\widehat{\textbf{n}}}}\right)$. Finally, we have shown that 
			$\sup_{y\in \mathbb{R}}\Vert\widehat{M}_{\textbf{n}}(y)-M(y)\Vert=O_p\left(b_{\textbf{n}}^{k}+  \dfrac{1}{b_{\textbf{n}}}\sqrt{\dfrac{\log \widehat{\textbf{n}}}{\widehat{\textbf{n}}}}\right) $.		
	
	\subsection{Proof of Theorem 3.1}
		We have 
		$\widehat\Gamma_{\textbf{n}}-\Gamma=2(\widehat\Sigma_{\textbf{n}}-\Sigma)+(\widehat\Lambda_{\textbf{n}}-\Lambda)
		$.
	Loubes and Yao (2013), show that $\widehat\Sigma_{\textbf{n}}-\Sigma=O_p\left(b_{\textbf{n}}^k+\dfrac{\phi_{\textbf{n}}^2}{e_{\textbf{n}}^2} \right)$; it remains to treat the second term. Putting 
	\[
f_{e,\textbf{n}}(y)=\max\{f(y),e_{\textbf{n}}\},\,\,
	R_{e,\textbf{n}}(y)=\dfrac{M(y)}{f_{e,\textbf{n}}(y)},\,\,
	r_{e,\textbf{n}}(y)=\dfrac{m(y)}{f_{e,\textbf{n}}(y)}\,\,
\]
and
\[
C_{e,\textbf{n}}(y)=R_{e,\textbf{n}}(y)-r_{e,\textbf{n}}(y)r_{e,\textbf{n}}(y)^T,
\]
we have
\[
\widehat{\Lambda}_{\textbf{n}}-\Lambda=\dfrac{1}{\widehat{\textbf{n}}}\sum\limits_{\textbf{i}\in \mathcal{I}_{\textbf{n}}}\widehat{C}_{\textbf{n}}(Y_{\textbf{i}})^2-\mathbb{E}\left[ Cov(Z|Y)^2\right]=A_{1\textbf{n}}+A_{2\textbf{n}}+A_{3\textbf{n}},
	\]
where 
\[
A_{1\textbf{n}}=\dfrac{1}{\widehat{\textbf{n}}}\sum\limits_{\textbf{i}\in \mathcal{I}_{\textbf{n}}}C(Y_{\textbf{i}})^2-\mathbb{E}\left[ Cov(Z|Y)^2\right],\,\,
A_{2\textbf{n}}=\dfrac{1}{\widehat{\textbf{n}}}\sum\limits_{\textbf{i}\in \mathcal{I}_{\textbf{n}}}C_{e,\textbf{n}}(Y_{\textbf{i}})^2-\dfrac{1}{\widehat{\textbf{n}}}\sum\limits_{\textbf{i}\in \mathcal{I}_{\textbf{n}}}C(Y_{\textbf{i}})^2
\]
and
\[
A_{3\textbf{n}}=\dfrac{1}{\widehat{\textbf{n}}}\sum\limits_{\textbf{i}\in \mathcal{I}_{\textbf{n}}}\widehat{C}_{\textbf{n}}(Y_{\textbf{i}})^2-\dfrac{1}{\widehat{\textbf{n}}}\sum\limits_{\textbf{i}\in \mathcal{I}_{\textbf{n}}}C_{e,\textbf{n}}(Y_{\textbf{i}})^2.
\]

\noindent $\bullet$ Control on  $A_{1\textbf{n}}$

\bigskip

\noindent Using  Assumption 3.5, we have
\begin{eqnarray*}
\mathbb{E}\left(\Vert C(Y_{\textbf{i}})^2-\mathbb{E}\left[ Cov(Z|Y)^2\right]\Vert^{2+\theta}\right)&\leq&\mathbb{E}\left(\left(\Vert C(Y_{\textbf{i}})\Vert^2+\Vert C(Y)\Vert^2\right)^{2+\theta}\right)\\
&\leq&\mathbb{E}\bigg[\bigg(\left(\mathbb{E}\left(\Vert Z\Vert^2\vert Y_{\textbf{i}}\right)+\mathbb{E}\left(\Vert Z\Vert\vert Y_{\textbf{i}}\right)^2\right)^2\\
& &+\left(\mathbb{E}\left(\Vert Z\Vert^2\vert Y\right)+\mathbb{E}\left(\Vert Z\Vert\vert Y\right)^2\right)^2\bigg)^{2+\theta}\bigg]\\
&\leq& \left( 8D^4\right)^{2+\theta}<+\infty,
\end{eqnarray*}
and since it is a continuous function of  $W_{\textbf{i}}$, the process $(C(Y_{\textbf{i}})^2-\mathbb{E}(C(Y)^2))_{\textbf{i}}$ is also strongly mixing. 
In addition,  since  $\theta>2L\geq 2$ we have the inequality $\theta^2/(2+\theta)>1$ and, consequently, we obtain that $\sum \alpha(\widehat{\textbf{n}})^{\theta/(2+\theta)}\leq C\sum \widehat{\textbf{n}}^{-\theta^2/(2+\theta)}<+\infty$. Using Lemma 1, we conclude  that $A_{1\textbf{n}}=O_p(1/\widehat{\textbf{n}})=O_p\left( \dfrac{1}{\sqrt{\widehat{\textbf{n}}}}\right) $.

\bigskip

\noindent $\bullet$ Control on  $A_{2\textbf{n}}$

\bigskip
		
\noindent We have		 
			\begin{align*}
			A_{2\textbf{n}}&=\dfrac{1}{\widehat{\textbf{n}}}\sum\limits_{\textbf{i}\in \mathcal{I}_{\textbf{n}}}[C(Y_{\textbf{i}})-C_{e,\textbf{n}}(Y_{\textbf{i}})]C(Y_{\textbf{i}})-\dfrac{1}{\widehat{\textbf{n}}}\sum\limits_{\textbf{i}\in \mathcal{I}_{\textbf{n}}}C_{e,\textbf{n}}(Y_{\textbf{i}})[C_{e,\textbf{n}}(Y_{\textbf{i}})-C(Y_{\textbf{i}})]=A_{21\textbf{n}}-A_{22\textbf{n}}.
			\end{align*}
The first term is
		 
			\begin{align*}
A_{21\textbf{n}}&=\dfrac{1}{\widehat{\textbf{n}}}\sum\limits_{\textbf{i}\in \mathcal{I}_{\textbf{n}}}[C(Y_{\textbf{i}})-C_{e,\textbf{n}}(Y_{\textbf{i}})]C(Y_{\textbf{i}})\\
			 &=\dfrac{1}{\widehat{\textbf{n}}}\sum\limits_{\textbf{i}\in \mathcal{I}_{\textbf{n}}}[R(Y_{\textbf{i}})-R_{e,\textbf{n}}(Y_{\textbf{i}})]C(Y_{\textbf{i}})+\dfrac{1}{\widehat{\textbf{n}}}\sum\limits_{\textbf{i}\in \mathcal{I}_{\textbf{n}}}[r_{e,\textbf{n}}(Y_{\textbf{i}})r_{e,\textbf{n}}^T(Y_{\textbf{i}})-r(Y_{\textbf{i}})r^T(Y_{\textbf{i}})]C(Y_{\textbf{i}})\\
			&=A_{211\textbf{n}}+A_{212\textbf{n}}.
			\end{align*}
			Since
	\begin{align*}
	\Vert A_{211\textbf{n}}\Vert &\leq\dfrac{1}{\widehat{\textbf{n}}}\sum\limits_{\textbf{i}\in \mathcal{I}_{\textbf{n}}}\Vert R(Y_{\textbf{i}})-R_{e,\textbf{n}}(Y_{\textbf{i}})\Vert\,\, \Vert C(Y_{\textbf{i}})\Vert\\
	&\leq \dfrac{1}{\widehat{\textbf{n}}}\sum\limits_{\textbf{i}\in \mathcal{I}_{\textbf{n}}}\Vert M(Y_{\textbf{i}})\Vert\,\, \Vert C(Y_{\textbf{i}})\Vert\,\, \bigg\vert\dfrac{1}{f(Y_{\textbf{i}})}-\dfrac{1}{f_{e,\textbf{n}}(Y_{\textbf{i}})}\bigg\vert\\
	&\leq \dfrac{1}{\widehat{\textbf{n}}}\sum\limits_{\textbf{i}\in \mathcal{I}_{\textbf{n}}}\Vert R(Y_{\textbf{i}})\Vert\,\, \Vert C(Y_{\textbf{i}})\Vert\,\,\textbf{1}_{\{f(Y_{\textbf{i}})<e_{\textbf{n}}\}}\\
	&\leq \dfrac{1}{\widehat{\textbf{n}}}\sum\limits_{\textbf{i}\in \mathcal{I}_{\textbf{n}}}\bigg(\Vert R(Y_{\textbf{i}})\Vert^2+\Vert R(Y_{\textbf{i}})\Vert\,\,\Vert r(Y_{\textbf{i}})\Vert^2\bigg)\,\,\textbf{1}_{\{f(Y_{\textbf{i}})<e_{\textbf{n}}\}}
\end{align*}
it follows
\[
	\mathbb{E}\left(\sqrt{\widehat{\textbf{n}}}\Vert A_{211\textbf{n}}\Vert\right)\leq \sqrt{\widehat{\textbf{n}}}\mathbb{E}\left(\Vert R(Y)\Vert^2\textbf{1}_{\{f(Y)<e_{\textbf{n}}\}}\right) +\sqrt{\widehat{\textbf{n}}} \mathbb{E}\left(\Vert R(Y)\Vert\,\,\Vert r(Y)\Vert^2\textbf{1}_{\{f(Y)<e_{\textbf{n}}\}}\right) 
	\]
	and, from Assumption 3.4 and Markov inequality we deduce that 
	$A_{211\textbf{n}}=O_p\left( 1/\sqrt{\widehat{\textbf{n}}}\right) $. Furthermore,  
	\begin{align*}
	\Vert A_{212\textbf{n}}\Vert&\leq\dfrac{1}{\widehat{\textbf{n}}}\sum\limits_{\textbf{i}\in \mathcal{I}_{\textbf{n}}}\Vert m(Y_{\textbf{i}})m^T(Y_{\textbf{i}})\Vert \,\, \Vert C(Y_{\textbf{i}})\Vert\,\, \bigg\vert\dfrac{1}{f^2(Y_{\textbf{i}})}-\dfrac{1}{f_{e,\textbf{n}}^2(Y_{\textbf{i}})}\bigg\vert
\end{align*}
and since
\[
\bigg\vert\dfrac{1}{f^2(Y_{\textbf{i}})}-\dfrac{1}{f_{e,\textbf{n}}^2(Y_{\textbf{i}})}\bigg\vert\leq \dfrac{1}{f^2(Y_{\textbf{i}})}\,\,\textbf{1}_{\{f(Y_{\textbf{i}})<e_{\textbf{n}}\}}
\]
it follows
\begin{align*}
	\Vert A_{212\textbf{n}}\Vert&\leq \dfrac{1}{\widehat{\textbf{n}}}\sum\limits_{\textbf{i}\in \mathcal{I}_{\textbf{n}}}\Vert r(Y_{\textbf{i}})r^T(Y_{\textbf{i}})\Vert \,\, \Vert C(Y_{\textbf{i}})\Vert\,\,\textbf{1}_{\{f(Y_{\textbf{i}})<e_{\textbf{n}}\}}\\
&\leq \dfrac{1}{\widehat{\textbf{n}}}\sum\limits_{\textbf{i}\in \mathcal{I}_{\textbf{n}}}\left(\Vert R(Y_{\textbf{i}})\Vert\,\, \Vert r(Y_{\textbf{i}})\Vert^2+\Vert r(Y_{\textbf{i}})\Vert^4 \right)\,\,\textbf{1}_{\{f(Y_{\textbf{i}})<e_{\textbf{n}}\}}.
\end{align*}
Hence
\begin{align*}
	\mathbb{E}\left( \sqrt{\widehat{\textbf{n}}}\Vert A_{212\textbf{n}}\Vert\right)&\leq \sqrt{\widehat{\textbf{n}}}\mathbb{E}\left(\Vert r(Y)\Vert^4\textbf{1}_{\{f(Y)<e_{\textbf{n}}\}}\right) +\sqrt{\widehat{\textbf{n}}} \mathbb{E}\left(\Vert R(Y)\Vert\,\,\Vert r(Y)\Vert^2\textbf{1}_{\{f(Y)<e_{\textbf{n}}\}}\right) 
	\end{align*}
and, from Assumption 3.4 and Markov inequality, we deduce that  $A_{212\textbf{n}}=O_p\left( 1/\sqrt{\widehat{\textbf{n}}} \right)$.
	Then, we have $A_{21\textbf{n}}=O_p\left( 1/\sqrt{\widehat{\textbf{n}}} \right)$. On the other hand,
	\begin{align*}
	\Vert  A_{22\textbf{n}}\Vert&\leq \bigg\Vert  \dfrac{1}{\widehat{\textbf{n}}}\sum\limits_{\textbf{i}\in \mathcal{I}_{\textbf{n}}}C_{e,\textbf{n}}(Y_{\textbf{i}})[C_{e,\textbf{n}}(Y_{\textbf{i}})-C(Y_{\textbf{i}})]\bigg\Vert \\
	 &\leq  \dfrac{1}{\widehat{\textbf{n}}}\sum\limits_{\textbf{i}\in \mathcal{I}_{\textbf{n}}}\Vert C_{e,\textbf{n}}(Y_{\textbf{i}})\Vert\,\, \Vert R_{e,\textbf{n}}(Y_{\textbf{i}})-R(Y_{\textbf{i}})\Vert+  \dfrac{1}{\widehat{\textbf{n}}}\sum\limits_{\textbf{i}\in \mathcal{I}_{\textbf{n}}}\Vert C_{e,\textbf{n}}(Y_{\textbf{i}})\Vert\,\, \Vert r_{e,\textbf{n}}(Y_{\textbf{i}})r_{e,\textbf{n}}^T(Y_{\textbf{i}})-r(Y_{\textbf{i}})r^T(Y_{\textbf{i}})\Vert.
\end{align*}
Since
\[
\Vert R_{e,\textbf{n}}(Y_{\textbf{i}})-R(Y_{\textbf{i}})\Vert=\Vert M(Y_{\textbf{i}})\Vert\, \bigg\vert\dfrac{1}{f(Y_{\textbf{i}})}-\dfrac{1}{f_{e,\textbf{n}}(Y_{\textbf{i}})}\bigg\vert
\leq \Vert R(Y_{\textbf{i}})\Vert\,\textbf{1}_{\{f(Y)<e_{\textbf{n}}\}},
\]
\[
\Vert C_{e,\textbf{n}}(Y_{\textbf{i}})\Vert\leq  \Vert R_{e,\textbf{n}}(Y_{\textbf{i}})\Vert+\Vert r_{e,\textbf{n}}(Y_{\textbf{i}})\Vert^2\leq \Vert R(Y_{\textbf{i}})\Vert+\Vert r(Y_{\textbf{i}})\Vert^2
\]
and
\begin{eqnarray*}
 \Vert r_{e,\textbf{n}}(Y_{\textbf{i}})r_{e,\textbf{n}}^T(Y_{\textbf{i}})-r(Y_{\textbf{i}})r^T(Y_{\textbf{i}})\Vert&=&\Vert  m(Y_{\textbf{i}})\, m(Y_{\textbf{i}})^T\Vert\,\,
\bigg\vert\dfrac{1}{f^2(Y_{\textbf{i}})}-\dfrac{1}{f_{e,\textbf{n}}^2(Y_{\textbf{i}})}\bigg\vert\\
&\leq&  \dfrac{\Vert  m(Y_{\textbf{i}})\, m(Y_{\textbf{i}})^T\Vert}{f^2(Y_{\textbf{i}})}\,\,\textbf{1}_{\{f(Y_{\textbf{i}})<e_{\textbf{n}}\}}\\
&=&\Vert r(Y_{\textbf{i}})\Vert^2\,\,\textbf{1}_{\{f(Y_{\textbf{i}})<e_{\textbf{n}}\}},
\end{eqnarray*}
it follows
\begin{align*}	
	\mathbb{E}&\left( \sqrt{\widehat{\textbf{n}}}\Vert A_{22\textbf{n}}\Vert\right)\leq  \sqrt{\widehat{\textbf{n}}}\mathbb{E}\left(\Vert R(Y)\Vert^2\textbf{1}_{\{f(Y)<e_{\textbf{n}}\}}\right) +2\sqrt{\widehat{\textbf{n}}} \mathbb{E}\left(\Vert R(Y)\Vert\,\,\Vert r(Y)\Vert^2\textbf{1}_{\{f(Y)<e_{\textbf{n}}\}}\right)\\
&+
 \sqrt{\widehat{\textbf{n}}} \mathbb{E}\left(\Vert r(Y)\Vert^4\textbf{1}_{\{f(Y)<e_{\textbf{n}}\}}\right).
\end{align*}
 from Assumption 3.4 and Markov inequality we deduce that $A_{22\textbf{n}}=O_p\left( 1/\sqrt{\widehat{\textbf{n}}} \right)$, and we can conclude that $A_{2\textbf{n}}=O_p\left( 1/\sqrt{\widehat{\textbf{n}}} \right)$.

\bigskip

\noindent $\bullet$ Control on  $A_{3\textbf{n}}$

\bigskip
	
		\begin{align*}
		A_{3\textbf{n}}&=\dfrac{1}{\widehat{\textbf{n}}}\sum\limits_{\textbf{i}\in \mathcal{I}_{\textbf{n}}}[C_{e,\textbf{n}}(Y_{\textbf{i}})-\widehat C_{e,\textbf{n}}(Y_{\textbf{i}})]C_{e,\textbf{n}}(Y_{\textbf{i}})-\dfrac{1}{\widehat{\textbf{n}}}\sum\limits_{\textbf{i}\in \mathcal{I}_{\textbf{n}}}[\widehat C_{e,\textbf{n}}(Y_{\textbf{i}})-C_{e,\textbf{n}}(Y_{\textbf{i}})][\widehat C_{e,\textbf{n}}(Y_{\textbf{i}})-C_{e,\textbf{n}}(Y_{\textbf{i}})]+\\
		&\dfrac{1}{\widehat{\textbf{n}}}\sum\limits_{\textbf{i}\in \mathcal{I}_{\textbf{n}}}C_{e,\textbf{n}}(Y_{\textbf{i}})[C_{e,\textbf{n}}(Y_{\textbf{i}})-\widehat C_{e,\textbf{n}}(Y_{\textbf{i}})]=A_{31\textbf{n}}-A_{32\textbf{n}}+A_{33\textbf{n}}.
		\end{align*}
First,
	\begin{align*}
	A_{31\textbf{n}}&=\dfrac{1}{\widehat{\textbf{n}}}\sum\limits_{\textbf{i}\in \mathcal{I}_{\textbf{n}}}[R_{e,\textbf{n}}(Y_{\textbf{i}})-\widehat R_{e,\textbf{n}}(Y_{\textbf{i}})]C_{e,\textbf{n}}(Y_{\textbf{i}})+\dfrac{1}{\widehat{\textbf{n}}}\sum\limits_{\textbf{i}\in \mathcal{I}_{\textbf{n}}}[\widehat r_{e,\textbf{n}}(Y_{\textbf{i}})\widehat r_{e,\textbf{n}}^T(Y_{\textbf{i}})-r_{e,\textbf{n}}(Y_{\textbf{i}})r_{e,\textbf{n}}^T(Y_{\textbf{i}})]C_{e,\textbf{n}}(Y_{\textbf{i}})\\
	&=\dfrac{1}{\widehat{\textbf{n}}}\sum\limits_{\textbf{i}\in \mathcal{I}_{\textbf{n}}}[R_{e,\textbf{n}}(Y_{\textbf{i}})-\widehat R_{e,\textbf{n}}(Y_{\textbf{i}})]C_{e,\textbf{n}}(Y_{\textbf{i}})+\dfrac{1}{\widehat{\textbf{n}}}\sum\limits_{\textbf{i}\in \mathcal{I}_{\textbf{n}}}[\widehat r_{e,\textbf{n}}(Y_{\textbf{i}})-r_{e,\textbf{n}}(Y_{\textbf{i}})][\widehat r_{e,\textbf{n}}(Y_{\textbf{i}})-r_{e,\textbf{n}}(Y_{\textbf{i}})]^TC_{e,\textbf{n}}(Y_{\textbf{i}})\\
	&+\dfrac{1}{\widehat{\textbf{n}}}\sum\limits_{\textbf{i}\in \mathcal{I}_{\textbf{n}}}[\widehat r_{e,\textbf{n}}(Y_{\textbf{i}})-r_{e,\textbf{n}}(Y_{\textbf{i}})] r_{e,\textbf{n}}(Y_{\textbf{i}})^TC_{e,\textbf{n}}(Y_{\textbf{i}})+\dfrac{1}{\widehat{\textbf{n}}}\sum\limits_{\textbf{i}\in \mathcal{I}_{\textbf{n}}} r_{e,\textbf{n}}(Y_{\textbf{i}})[\widehat r_{e,\textbf{n}}(Y_{\textbf{i}})-r_{e,\textbf{n}}(Y_{\textbf{i}})]^TC_{e,\textbf{n}}(Y_{\textbf{i}})\\
	&=A_{311\textbf{n}}+A_{312\textbf{n}}+A_{313\textbf{n}}+A_{314\textbf{n}}
	\end{align*}
		and
	\begin{align*}
	\Vert A_{311\textbf{n}}\Vert&\leq \dfrac{1}{\widehat{\textbf{n}}}\sum\limits_{\textbf{i}\in \mathcal{I}_{\textbf{n}}}\Vert R_{e,\textbf{n}}(Y_{\textbf{i}})-\widehat R_{e,\textbf{n}}(Y_{\textbf{i}})\Vert \,\,\Vert C_{e,\textbf{n}}(Y_{\textbf{i}})\Vert\\
	&\leq \dfrac{1}{\widehat{\textbf{n}}}\sum\limits_{\textbf{i}\in \mathcal{I}_{\textbf{n}}}\bigg\Vert  \dfrac{{R}_{e,\textbf{n}}(Y_{\textbf{i}})}{\widehat{f}_{e,\textbf{n}}(Y_{\textbf{i}})}\left( f_{e,\textbf{n}}(Y_{\textbf{i}})-\widehat{f}_{e,\textbf{n}}(Y_{\textbf{i}})\right) +\dfrac{1}{\widehat{f}_{e,\textbf{n}}(Y_{\textbf{i}})}\left( \widehat{M}_{\textbf{n}}(Y_{\textbf{i}})-M(Y_{\textbf{i}})\right) \bigg\Vert\,\,  \Vert C_{e,\textbf{n}}(Y_{\textbf{i}})\Vert.
\end{align*}
We have
$
\dfrac{1}{\widehat{f}_{e,\textbf{n}}(y)}\leq \dfrac{1}{e_{\textbf{n}}}
$
and Assumption 3.5 implies $\Vert C_{e,\textbf{n}}(Y_{\textbf{i}})\Vert\leq 2D^2$  and $\Vert R(Y_{\textbf{i}})\Vert\leq D^2$. Furthermore, we have from Zhu and Fang (1996):  $\Vert f_{e,\textbf{n}}-\widehat{f}_{e,\textbf{n}} \Vert_{\infty}\leq \Vert f-\widehat{f}_{\textbf{n}} \Vert_{\infty}$. Thus,
\begin{align*}
	\Vert A_{311\textbf{n}}\Vert&\leq \Vert f-\widehat{f}_{\textbf{n}} \Vert_{\infty}\dfrac{1}{\widehat{\textbf{n}}e_{\textbf{n}}}\sum\limits_{\textbf{i}\in \mathcal{I}_{\textbf{n}}}\Vert R(Y_{\textbf{i}})\Vert\,\, \Vert C_{e,\textbf{n}}(Y_{\textbf{i}})\Vert+ \Vert M-\widehat{M}_{\textbf{n}} \Vert_{\infty}\dfrac{1}{\widehat{\textbf{n}}e_{\textbf{n}}}\sum\limits_{\textbf{i}\in \mathcal{I}_{\textbf{n}}} \Vert C_{e,\textbf{n}}(Y_{\textbf{i}})\Vert \\
&\leq \dfrac{2D^4}{e_{\textbf{n}}}\,\Vert f-\widehat{f}_{\textbf{n}} \Vert_{\infty}+ \dfrac{2D^2}{e_{\textbf{n}}} \Vert M-\widehat{M}_{\textbf{n}} \Vert_{\infty},
	\end{align*}
	and using Lemma 3 and Lemma 4, we obtain $A_{311\textbf{n}}=O_p\left( \dfrac{\phi_{\textbf{n}}}{e_{\textbf{n}}} \right)$.
	Similar argument leads  us to get
	$A_{313\textbf{n}}=O_p\left( \dfrac{\phi_{\textbf{n}}}{e_{\textbf{n}}} \right)$ and $A_{314\textbf{n}}=O_p\left( \dfrac{\phi_{\textbf{n}}}{e_{\textbf{n}}} \right)$. On the other hand,
	\begin{align*}
	\Vert A_{312\textbf{n}}\Vert&\leq\dfrac{1}{\widehat{\textbf{n}}}\sum\limits_{\textbf{i}\in \mathcal{I}_{\textbf{n}}}\bigg\Vert  \dfrac{r_{e,\textbf{n}}(Y_{\textbf{i}})}{\widehat{f}_{e,\textbf{n}}(Y_{\textbf{i}})}\left( f_{e,\textbf{n}}(Y_{\textbf{i}})-\widehat{f}_{e,\textbf{n}}(Y_{\textbf{i}})\right) -\dfrac{1}{\widehat{f}_{e,\textbf{n}}(Y_{\textbf{i}})}\left( \widehat{m}_{\textbf{n}}(Y_{\textbf{i}})-m(Y_{\textbf{i}})\right) \bigg\Vert ^2\,\, \Vert C_{e,\textbf{n}}(Y_{\textbf{i}})\Vert\\
	&= \dfrac{1}{\widehat{\textbf{n}}}\sum\limits_{\textbf{i}\in \mathcal{I}_{\textbf{n}}}\bigg\Vert   \dfrac{r_{e,\textbf{n}}(Y_{\textbf{i}})}{\widehat{f}_{e,\textbf{n}}(Y_{\textbf{i}})}\left( f_{e,\textbf{n}}(Y_{\textbf{i}})-\widehat{f}_{e,\textbf{n}}(Y_{\textbf{i}})\right)\bigg\Vert ^2\Vert C_{e,\textbf{n}}(Y_{\textbf{i}})\Vert  +\dfrac{1}{\widehat{\textbf{n}}}\sum\limits_{\textbf{i}\in \mathcal{I}_{\textbf{n}}}\bigg\Vert \dfrac{\widehat{m}_{\textbf{n}}(Y_{\textbf{i}})-m(Y_{\textbf{i}})}{\widehat{f}_{e,\textbf{n}}(Y_{\textbf{i}})}\bigg\Vert ^2  \Vert C_{e,\textbf{n}}(Y_{\textbf{i}})\Vert\\
	&-\dfrac{2}{\widehat{\textbf{n}}}\sum\limits_{\textbf{i}\in \mathcal{I}_{\textbf{n}}}\Vert C_{e,\textbf{n}}(Y_{\textbf{i}})\Vert \dfrac{ f_{e,\textbf{n}}(Y_{\textbf{i}})-\widehat{f}_{e,\textbf{n}}(Y_{\textbf{i}})}{\widehat{f}_{e,\textbf{n}}^2(Y_{\textbf{i}})} \,\,r_{e,\textbf{n}}(Y_{\textbf{i}})^T\left( \widehat{m}_{\textbf{n}}(Y_{\textbf{i}})-m(Y_{\textbf{i}})\right);
	\end{align*}
	from Zhu and Fang(1996), we have that $\vert f_{e,\textbf{n}}(Y_{\textbf{i}})-\widehat{f}_{e,\textbf{n}}(Y_{\textbf{i}})\vert\leq \vert \widehat{f}_{\textbf{n}}(Y_{\textbf{i}})-f(Y_{\textbf{i}})\vert\leq \Vert\widehat{f}_{\textbf{n}}-f\Vert_{\infty}$, then
	\begin{align*}
	\Vert A_{312\textbf{n}}\Vert&\leq \dfrac{\Vert f-\widehat{f}_{\textbf{n}} \Vert_{\infty}^2}{e_{\textbf{n}}^2}\left( \dfrac{1}{\widehat{\textbf{n}}}\sum\limits_{\textbf{i}\in \mathcal{I}_{\textbf{n}}}\Vert r(Y_{\textbf{i}})\Vert^2\,\,\Vert C_{e,\textbf{n}}(Y_{\textbf{i}})\Vert\right)+\dfrac{\Vert m-\widehat{m}_{\textbf{n}} \Vert_{\infty}^2}{e_{\textbf{n}}^2}\left( \dfrac{1}{\widehat{\textbf{n}}}\sum\limits_{\textbf{i}\in \mathcal{I}_{\textbf{n}}}\Vert C_{e,\textbf{n}}(Y_{\textbf{i}})\Vert\right)  \\
	&+\dfrac{\Vert m-\widehat{m}_{\textbf{n}} \Vert_{\infty}\,\,\Vert f-\widehat{f}_{\textbf{n}} \Vert_{\infty} }{e_{\textbf{n}}^2}\left( \dfrac{1}{\widehat{\textbf{n}}}\sum\limits_{\textbf{i}\in \mathcal{I}_{\textbf{n}}}\Vert r(Y_{\textbf{i}})\Vert\,\, \Vert C_{e,\textbf{n}}(Y_{\textbf{i}})\Vert\right)\\
&\leq 2D^4\dfrac{\Vert f-\widehat{f}_{\textbf{n}} \Vert_{\infty}^2}{e_{\textbf{n}}^2}+2D^2\dfrac{\Vert m-\widehat{m}_{\textbf{n}} \Vert_{\infty}^2}{e_{\textbf{n}}^2}  +2D^3\,\dfrac{\Vert m-\widehat{m}_{\textbf{n}} \Vert_{\infty}\,\,\Vert f-\widehat{f}_{\textbf{n}} \Vert_{\infty} }{e_{\textbf{n}}^2}.
\end{align*}
	Using Lemma 3, we obtain $\Vert A_{312\textbf{n}}\Vert=O_p(\frac{\phi_{\textbf{n}}^2}{e_{\textbf{n}}^2})$ and, therefore, 
	$A_{31\textbf{n}}=O_p(\frac{\phi_{\textbf{n}}}{e_{\textbf{n}}})+O_p(\frac{\phi_{\textbf{n}}^2}{e_{\textbf{n}}^2})$. It can be noticed that $\Vert A_{33\textbf{n}}\Vert=\Vert A_{31\textbf{n}}\Vert$. Thus, 
	$A_{33\textbf{n}}=O_p(\frac{\phi_{\textbf{n}}}{e_{\textbf{n}}})+O_p(\frac{\phi_{\textbf{n}}^2}{e_{\textbf{n}}^2})$. Further, 
	\begin{align*}
	A_{32\textbf{n}}&=\dfrac{1}{\widehat{\textbf{n}}}\sum\limits_{\textbf{i}\in \mathcal{I}_{\textbf{n}}}[\widehat R_{e,\textbf{n}}(Y_{\textbf{i}})-R_{e,\textbf{n}}(Y_{\textbf{i}})][\widehat R_{e,\textbf{n}}(Y_{\textbf{i}})-R_{e,\textbf{n}}(Y_{\textbf{i}})]\\
	&-\dfrac{1}{\widehat{\textbf{n}}}\sum\limits_{\textbf{i}\in \mathcal{I}_{\textbf{n}}}[\widehat R_{e,\textbf{n}}(Y_{\textbf{i}})-R_{e,\textbf{n}}(Y_{\textbf{i}})][\widehat r_{e,\textbf{n}}(Y_{\textbf{i}})\widehat r_{e,\textbf{n}}^T(Y_{\textbf{i}})-r_{e,\textbf{n}}(Y_{\textbf{i}})r_{e,\textbf{n}}^T(Y_{\textbf{i}})]\\
	&-\dfrac{1}{\widehat{\textbf{n}}}\sum\limits_{\textbf{i}\in \mathcal{I}_{\textbf{n}}}[\widehat r_{e,\textbf{n}}(Y_{\textbf{i}})\widehat r_{e,\textbf{n}}^T(Y_{\textbf{i}})-r_{e,\textbf{n}}(Y_{\textbf{i}})r_{e,\textbf{n}}^T(Y_{\textbf{i}})][\widehat R_{e,\textbf{n}}(Y_{\textbf{i}})-R_{e,\textbf{n}}(Y_{\textbf{i}})]\\
	&+\dfrac{1}{\widehat{\textbf{n}}}\sum\limits_{\textbf{i}\in \mathcal{I}_{\textbf{n}}}[\widehat r_{e,\textbf{n}}(Y_{\textbf{i}})\widehat r_{e,\textbf{n}}^T(Y_{\textbf{i}})-r_{e,\textbf{n}}(Y_{\textbf{i}})r_{e,\textbf{n}}^T(Y_{\textbf{i}})]^2\\
	 &=	A_{321\textbf{n}}-A_{322\textbf{n}}-A_{323\textbf{n}}+	A_{324\textbf{n}}
	\end{align*}
		Similar arguments than those  used for   $A_{31\textbf{n}}$ lead to $A_{32\textbf{n}}=O_p(\frac{\phi_{\textbf{n}}}{e_{\textbf{n}}})+O_p(\frac{\phi_{\textbf{n}}^2}{e_{\textbf{n}}^2})+O_p(\frac{\phi_{\textbf{n}}^3}{e_{\textbf{n}}^3})+O_p(\frac{\phi_{\textbf{n}}^4}{e_{\textbf{n}}^4})$. We can then conclude that $A_{3\textbf{n}}=O_p(\frac{\phi_{\textbf{n}}}{e_{\textbf{n}}})+O_p(\frac{\phi_{\textbf{n}}^2}{e_{\textbf{n}}^2})+O_p(\frac{\phi_{\textbf{n}}^3}{e_{\textbf{n}}^3})+O_p(\frac{\phi_{\textbf{n}}^4}{e_{\textbf{n}}^4})$. Finally, the theorem is proven.
			\subsection{Proof of Corollary 3.2}
		Assumption 3.7 leads to $\dfrac{\widehat{\textbf{n}}^{1/2}}{e_{\textbf{n}}}\phi_{\textbf{n}}\sim   \widehat{\textbf{n}}^{-1/2+c_2+c_1}\sqrt{\log\widehat{\textbf{n}}}$.  Since  $-1/2+c_2+c_1<-1/2+2c_2+c_1<0$, we obtain that  $\dfrac{\widehat{\textbf{n}}^{1/2}}{e_{\textbf{n}}}\phi_{\textbf{n}}\rightarrow 0$ as  $\textbf{n}\rightarrow +\infty$. Then, $\frac{\phi_{\textbf{n}}}{e_{\textbf{n}}}=O_p(\widehat{\textbf{n}}^{-1/2})=o_p(1)$ and, consequently, $\frac{\phi_{\textbf{n}}^\ell}{e_{\textbf{n}}^\ell}=O_p\left( \frac{\phi_{\textbf{n}}}{e_{\textbf{n}}}\right)=O_p(\widehat{\textbf{n}}^{-1/2})$ for $\ell =2,3,4$. This gives the result.
			\subsection{Proof of Corollary 3.3}
\noindent
Clearly,  
\[
\widehat{\Sigma}_{\textbf{n}}-\Sigma= \dfrac{1}{\widehat{\textbf{n}}}\sum\limits_{\textbf{i}\in \mathcal{I}_{\textbf{n}}}\mathcal{V}_{\textbf{i}}- \left(\overline{X}-\mathbb{E}(X)\right)\overline{X}^T-\mathbb{E}(X)\left(\overline{X}-\mathbb{E}(X)\right)^T,
\]
where   $\mathcal{V}_{\textbf{i}}:=X_{\textbf{i}}X_{\textbf{i}}^T- \mathbb{E}(X) \mathbb{E}(X)^T$. Using  Assumption 3.5, we have
\begin{eqnarray*}
\mathbb{E}\left(\Vert \mathcal{V}_{\textbf{i}}\Vert^{2+\theta}\right)&\leq&\mathbb{E}\left(\left(\Vert  X_{\textbf{i}}\Vert^2+\Vert\mathbb{E}(X)\Vert^2\right)^{2+\theta}\right)\\
&\leq&\mathbb{E}\left( \left(\left\Vert \Sigma\Vert_\infty\,\Vert Z_{\textbf{i}}\Vert+\mathbb{E}\left(\Vert X\Vert^2\right)\right)^2+\mathbb{E}\left(\Vert X\Vert^2\right)\right)^{2+\theta}\right)\\
&\leq& \left( 2\Vert \Sigma\Vert_\infty^2\,D^2+3\mathbb{E}\left(\Vert X\Vert^2\right)\right)^{2+\theta}<+\infty,
\end{eqnarray*}
and
\begin{eqnarray*}
\mathbb{E}\left(\Vert X_{\textbf{i}}-\mathbb{E}(X)\Vert^{2+\theta}\right)&\leq&\mathbb{E}\left(\left(\Vert  X_{\textbf{i}}\Vert+\mathbb{E}(\Vert X\Vert)\right)^{2+\theta}\right)\\
&\leq&\mathbb{E}\left(\left( \Vert \Sigma\Vert_\infty\,\Vert Z_{\textbf{i}}\Vert+2\mathbb{E}\left(\Vert X\Vert\right)\right)^{2+\theta}\right)\\
&\leq& \left( \Vert \Sigma\Vert_\infty\,D+2\mathbb{E}\left(\Vert X\Vert\right)\right)^{2+\theta}<+\infty.
\end{eqnarray*}
In addition,  since  $\theta>2L\geq 2$ we have the inequality $\theta^2/(2+\theta)>1$ and, consequently, we obtain that $\sum \alpha(\widehat{\textbf{n}})^{\theta/(2+\theta)}\leq C\sum \widehat{\textbf{n}}^{-\theta^2/(2+\theta)}<+\infty$. Using Lemma 1, we obtain  $\dfrac{1}{\widehat{\textbf{n}}}\sum\limits_{\textbf{i}\in \mathcal{I}_{\textbf{n}}}\mathcal{V}_{\textbf{i}}=O_p(1/\widehat{\textbf{n}})=o_p(1)$ and $\overline{X}-\mathbb{E}(X)=O_p(1/\widehat{\textbf{n}})=o_p(1)$. Therefore,  $\widehat{\Sigma}_{\textbf{n}}-\Sigma=O_p(1/\widehat{\textbf{n}})=o_p(1)$. From the continuity of the map  $A\longmapsto A^{-1/2}$ we then deduce that $\Vert\widehat{\Sigma}_{\textbf{n}}^{-1/2}-\Sigma^{-1/2}\Vert=o_p(1)$.  Applying Lemma 1 of Ferr\'e  and Yao(2003)  gives  $\Vert \widehat{\tau}_{j}-\tau_j\Vert\leq b_j \Vert \widehat{\Gamma}_{\textbf{n}}-\Gamma \Vert_{\infty}$, where 
 $b_1=2\sqrt{2}/(\lambda_1-\lambda_2)$ and $b_j=2\sqrt{2}/\min\{\lambda_{j-1}-\lambda_{j};\lambda_{j}-\lambda_{j+1}\}$ for $j\geq 2$, we have for $j=1,2,\cdots, N$. Then,  Theorem 3.1 permits to conclude that $\Vert \widehat{\tau}_{j}-\tau_j\Vert =o_p(1) $  for $j=1,2,\cdots, N$. Finally, from
		$
		\widehat{\beta}_j-\beta_j=\left(\widehat{\Sigma}_{\textbf{n}}^{-1/2} -\Sigma^{-1/2}\right)\widehat{\tau}_j+\Sigma^{-1/2}\left(\widehat{\tau}_j-\tau_j\right)
$, 
we deduce that
\[
		\Vert\widehat{\beta}_j-\beta_j\Vert\leq\Vert\widehat{\Sigma}_{\textbf{n}}^{-1/2}-\Sigma^{-1/2}\Vert\,\, \Vert\widehat{\tau}_j\Vert+\Vert\Sigma^{-1/2}\Vert\,\, \Vert\widehat{\tau}_j-\tau_j\Vert
\]
and we conclude that 
		$\Vert\widehat{\beta}_j-\beta_j\Vert=o_p(1)$.

\end{document}